\documentclass[final,3p]{elsarticle}

\usepackage{amssymb}
\usepackage{amsmath}
\usepackage{xcolor}
\usepackage{lineno}

\usepackage{hyperref}
\hypersetup{hidelinks,colorlinks=false,breaklinks=true,urlcolor=primary,bookmarksopen=false}

\hypersetup{pdfinfo={Title={Hejlesen-2020}}}

 \biboptions{square,comma}

\journal{arXiv}

\begin{document}

\begin{frontmatter}

\title{Non-singular Green's functions for the unbounded Poisson equation in 
one, two and three dimensions}

\author[DTU]{Mads M{\o}lholm Hejlesen}
\author[UCL]{Gr{\'e}goire Winckelmans}
\author[DTU,ETH]{Jens Honor\'{e} Walther\corref{cor1}}

\address[DTU]{Department of Mechanical Engineering, 
Technical University of Denmark, Kgs.\ Lyngby, Denmark}
\address[UCL]{Institute of Mechanics, Materials and Civil Engineering, 
Universit\'{e} catholique de Louvain, Louvain-la-Neuve, Belgium}
\address[ETH]{Computational Science and Engineering Laboratory, 
ETH Z\"{u}rich, Z\"{u}rich, Switzerland}

\cortext[cor1]{Corresponding author at: 
Department of Mechanical Engineering, Technical University of Denmark,
Building 403, DK-2800 Kgs.\ Lyngby, Denmark. 
Tel.: + 45 4525 4327; fax: + 45 4588 4325. 
E-mail address: jhw@mek.dtu.dk (J.~H.~Walther).}

\begin{abstract}

This paper is a revised version of the original paper of same title---published 
in Applied Mathematics Letters 89 \cite{Hejlesen:2019}---containing some 
corrections and clarifications to the original text.

We derive non-singular Green's functions for the unbounded 
Poisson equation in one, two and three dimensions, using a 
cut-off function in the Fourier domain to impose a smallest 
length scale when deriving the Green's function. 
The resulting non-singular Green's functions are relevant to applications 
which are restricted to a minimum resolved length scale 
(e.g. a mesh size $h$) and thus cannot handle the singular Green's function 
of the continuous Poisson equation. 
We furthermore derive the gradient vector of the non-singular Green's 
function, as this is useful in applications where the 
Poisson equation represents potential functions of a vector field.  

\end{abstract}


\end{frontmatter}

\section{Introduction}

The use of Green's functions for solving 
linear differential equations in an unbounded domain, i.e. with free-space 
boundary conditions, is a frequently used methodology. 
The Green's function represents the impulse response function to a linear 
differential operator, such as the Laplace operator in the case of the 
Poisson equation. 
Once obtained, the Green's function can be used, by utilizing the 
superposition principle to obtain the solution of an inhomogeneous equation 
by convolving the Green's function with the right-hand-side field of the 
equation.

The Green's function of the Laplace operator is singular 
at its origin. Applied in discretized numerical calculations, 
the singularity of the Green's function evidently causes a number of difficulties. 
In order to amend this, smoothing regularization techniques have 
been applied (e.g.\ \cite{Leonard:1980,Beale:1985,Winckelmans:1993}) 
which imposes the effect of a continuous and smooth field distribution 
around the discrete points, 
which in effect can be used to avoid the singularity of the Green's function. 
However, most regularization methods that have been applied are based on 
functions that only conserve a finite number of field moments, and are thus 
only accurate up to a finite order of convergence rate.

Vico et al.\ \cite{Vico:2016} derived non-singular Green's function in the 
two and three dimensional Fourier domain by imposing an isotropic 
maximum length scale of the real domain. 
This method was shown to provide a spectral accuracy in computations, where 
the real domain was sufficiently extended in all directions to provide a 
sufficient resolution when evaluating the Green's function in 
the Fourier domain.
However, as the Green's function is a radial function the imposed maximum length 
scale dictates the resolution criteria in the Fourier domain, and thus the 
length of the integrated domain, for all directions.
This makes the method of Vico et al.\ \cite{Vico:2016} potentially 
inefficient for elongated domains, as the imposed length scale must be proportional 
to the maximum length of the domain to insure that the integration includes 
the full domain.

In this work we impose a minimum length scale in the real domain, 
analogous to a discretization length, in order to obtain a non-singular 
Green's functions in the real domain.
This is done using a cut-off function in the Fourier domain to impose 
a minimum length scale, allowing us to derive the Green's function analytically 
by Fourier analysis. 
As no explicit smoothing function is applied, the obtained non-singular 
Green's function gives an optimal accuracy when used in numerical calculations, 
as it is only subject to quadrature errors, converging in a rate proportional 
to the smoothness of the right-hand-side field of the Poisson equation.

As with the method of Vico et al.\ \cite{Vico:2016} the presented method also 
imposes an equivalent length scale for all directions, albeit here as a minimum 
length scale of the solution. 
Consequently, if the discretization is non-isotropic the proposed method 
will provide a solution, in all directions, which limited to the minimum 
length scale of the lowest resolved direction.

The gradient vector of the non-singular Green's function is also derived. 
This function is useful in applications where the Poisson equation represents 
a vector or scalar potential function of a vector field as it enables a direct 
solution of the vector field. 

The proposed non-singular Green's functions may be used directly in efficient 
computational methods such as the fast multipole method 
\cite{Carrier:1988,Barnes:1986} or the mesh based FFT 
solver \cite{Hockney:1970}.
The mesh based Poisson solver was recently shown capable of obtaining a high 
order convergence rate by using regularized Green's functions 
\cite{Hejlesen:2013,Hejlesen:2015d,Hejlesen:2016}. 
The mesh based solver was extended to handle domains with a combination of 
unbounded and periodic directions by 
Chatelain and Koumoutsakos \cite{Chatelain:2010} 
and for high order accuracy in Spietz et al. \cite{Spietz:2018}. 
Here it is utilized that the domain is intrinsically bounded in the 
periodic directions, yielding a straightforward Green's function 
in the Fourier domain, which may be combined with the Fourier transformed 
unbounded Green's function of reduced dimensionality for the 
unbounded directions.

Using the non-singular Green's functions presented here in the mesh 
based Poisson solver effectively results in a solution of spectral 
accuracy for both unbounded domains as well as a domains with a combination of 
unbounded and periodic directions. 
The reader is referred to the open source software 
\cite{Hejlesen:2018} for a numerical implementation of this.

\section{Methodology}

The Poisson equation in an unbounded domain is formally stated as:
\begin{equation}
 \boldsymbol{\nabla}^2 A(\boldsymbol{x}) = -B(\boldsymbol{x})
\qquad \text{where} \qquad 
 A(\boldsymbol{x}) \rightarrow 0 
\quad \text{for} \quad 
|\boldsymbol{x}| \rightarrow \infty
 \label{eq:Poisson}
\end{equation}
Here $B(\boldsymbol{x})$ is a known compactly supported field 
and $A(\boldsymbol{x})$ is the desired solution field.
In many applications such as astrophysics, electrodynamics and vortex dynamics, 
the vector field to be solved $\boldsymbol{v}(\boldsymbol{x})$ is described 
by potential functions using the Helmholtz decomposition:
\begin{equation}
 \boldsymbol{v}(\boldsymbol{x}) = \boldsymbol{\nabla} \times \boldsymbol{\psi}(\boldsymbol{x}) - \boldsymbol{\nabla} \phi(\boldsymbol{x})
 \qquad \mbox{where} \qquad \boldsymbol{\nabla} \cdot \boldsymbol{\psi}(\boldsymbol{x}) = 0 
 \label{eq:Helmholtz}
\end{equation}
The fundamental operations describing the pointwise conservation of the flux 
and circulation of the vector field $\boldsymbol{v}$ is the 
divergence $\vartheta(\boldsymbol{x}) = \boldsymbol{\nabla} \cdot \boldsymbol{v}(\boldsymbol{x})$ and 
the curl $\boldsymbol{\omega}(\boldsymbol{x}) = \boldsymbol{\nabla} \times \boldsymbol{v}(\boldsymbol{x})$, 
respectively. From Eq.~\eqref{eq:Helmholtz} it follows that these may be 
expressed by the potential functions as:
\begin{equation}
 \vartheta(\boldsymbol{x}) = - \boldsymbol{\nabla}^2 \phi(\boldsymbol{x})
\qquad \mbox{and} \qquad
 \boldsymbol{\omega}(\boldsymbol{x}) = - \boldsymbol{\nabla}^2 \boldsymbol{\psi}(\boldsymbol{x})
 \label{eq:potential Poisson}
\end{equation}
Both potentials may thus be obtained by solving a Poisson equation. 
By utilizing the linearity of the Poisson equation, 
and considering the equations transformed into the Fourier domain, 
we obtain the algebraic equations:
\begin{equation}
\begin{array}{l c l}
 - k^2  \widehat{\phi}(\boldsymbol{k}) = - \widehat{\vartheta}(\boldsymbol{k})
 &  \Leftrightarrow 
 & \widehat{\phi}(\boldsymbol{k})  = \widehat{G}(\boldsymbol{k}) \, \widehat{\vartheta}(\boldsymbol{k})
 \\[0.3cm]
 - k^2  \widehat{\boldsymbol{\psi}}(\boldsymbol{k}) = - \widehat{\boldsymbol{\omega}}(\boldsymbol{k})
 & \Leftrightarrow 
 & \widehat{\boldsymbol{\psi}}(\boldsymbol{k}) = \widehat{G}(\boldsymbol{k}) \, \widehat{\boldsymbol{\omega}}(\boldsymbol{k})
\end{array}
\quad \text{where} \quad 
 \widehat{G}(\boldsymbol{k}) = \frac{1}{k^2}
 \quad \text{for} \quad k > 0
\label{eq:Fourier Poisson}
\end{equation}
Here ~$\widehat{\cdot}$~ denotes the field generated by the Fourier transform, 
$\boldsymbol{k}$ is the angular wave-number vector of the 
Fourier domain, and $k = |\boldsymbol{k}|$ is its radial length.
$\widehat{G}(\boldsymbol{k})$ is the Green's function in the Fourier domain 
representation.
We may furthermore obtain the vector field of Eq.~\eqref{eq:Helmholtz} directly 
by incorporating the gradient operator into the Green's function as 
$\boldsymbol{K}(\boldsymbol{x}) =\boldsymbol{\nabla} G(\boldsymbol{x})$, 
after which we obtain:
\begin{equation}
\begin{array}{l c l}
 - k^2 \, \widehat{\boldsymbol{v}}_\vartheta(\boldsymbol{k}) = - \text{i} \, \boldsymbol{k} \, \widehat{\vartheta}(\boldsymbol{k})
 & \Leftrightarrow 
 & \widehat{\boldsymbol{v}}_\vartheta(\boldsymbol{k})  = \widehat{\boldsymbol{K}}(\boldsymbol{k}) \, \widehat{\vartheta}(\boldsymbol{k})
 \\[0.3cm]
 - k^2 \, \widehat{\boldsymbol{v}}_\omega(\boldsymbol{k}) = - \text{i} \, \boldsymbol{k} \times \widehat{\boldsymbol{\omega}}(\boldsymbol{k})
 & \Leftrightarrow 
 & \widehat{\boldsymbol{v}}_\omega(\boldsymbol{k}) = \widehat{\boldsymbol{K}}(\boldsymbol{k}) \times \widehat{\boldsymbol{\omega}}(\boldsymbol{k})
\end{array}
\quad \text{where} \quad 
 \widehat{\boldsymbol{K}}(\boldsymbol{k}) = \frac{\text{i} \, \boldsymbol{k}}{k^2}
 \quad \text{for} \quad k > 0
\label{eq:Fourier Poisson vec}
\end{equation}

As the Green's function in the Fourier domain (Eq.~\eqref{eq:Fourier Poisson}) 
does not have a compact support, an unbounded Fourier domain is needed 
to obtain the exact solution. 
However, for a discrete approximation, such as that used in numerical simulations, 
the discretization of a function in the real domain, 
imposes an upper bound $k_s$ of the function image in the Fourier domain.
This upper bound is determined by the Nyquist-Shannon sampling theorem 
as $k_s = \pi/h$ being the highest resolved angular wave-number corresponding 
to the discretization length $h$.
Thus, the exact Green's function of Eq.~\eqref{eq:Fourier Poisson} cannot not be 
represented in a discretized domain, and using it will introduce errors 
in the calculation resulting in a limited convergence rate 
of $\mathcal{O}(h^2)$ \cite{Rasmussen:2011c}.
In order to amend this and derive a Green's function which is bounded by a 
maximum wavenumber, we regularize the Green's function using a radial 
cut-off function in the Fourier domain:
\begin{equation}
 \widehat{\zeta}(k) = \left\{ \begin{array}{l l l}
     1    & \text{for} & 0 < k \leq \dfrac{\pi}{h} \\[0.2cm]
     0    & \text{for} & k > \dfrac{\pi}{h}
   \end{array} \right.
 \label{eq:zeta}
\end{equation}
This leads to a regularized Green's function and a corresponding equation 
for the impulse response function:
\begin{equation}
 \widehat{G}(k) = \frac{\widehat{\zeta}(k)}{k^2} 
\quad \Leftrightarrow \quad
 \boldsymbol{\nabla}^2 G(r) = -\zeta(r)
 \label{eq:greens poisson 2}
\end{equation}
where $r = |\boldsymbol{x}|$ is the radial length of the position vector. 
Using the radial cut-off function of Eq.~\eqref{eq:zeta}, it is seen that the 
Green's function of the non-regularized equation (Eq.~\eqref{eq:Fourier Poisson}) 
is unchanged for $0 < k \leq \pi/h$.

\section{Results}

By utilizing the radial symmetry of the Green's function, the 
regularization function in the real domain $\zeta(\rho)$ may be obtained by 
considering the $d$-dimensional radial Fourier transform \cite{Stein:1971} 
for $d \geq 1$:
\begin{equation}
 \widehat{\zeta}(s) = (2\pi \sigma^2)^{\frac{d}{2}} \int_0^\infty 
 \zeta(\rho) \
\dfrac{J_{\frac{d}{2} - 1}(s \rho)}{(s \rho)^{\frac{d}{2} - 1}} \rho^{d - 1} \ d\rho
\label{eq:zetaFT}
\end{equation}
where we have introduced the normalized coordinates $s = \sigma k$ 
and $\rho = r/\sigma$ relative to the length scale $\sigma$.
The function $J_\nu$ is the Bessel function of the first kind and of 
order $\nu$. 
Using an integral property of the Bessel functions of the first kind 
\cite{Olver:2010}, 
we may satisfy Eq.~\eqref{eq:zeta} by: 
\begin{equation}
 \widehat{\zeta}(s) =
 \displaystyle\int_0^\infty \frac{ J_{\frac{d}{2}}(\rho) J_{\frac{d}{2}-1}(s \rho) }{s^{\frac{d}{2}-1}} \, d\rho 
 = \left\{ \begin{array}{l l l}
         1    & \text{for} & 0 < s < 1 \\[0.2cm]
 \frac{1}{2}  & \text{for} & s = 1 \\[0.2cm]
         0    & \text{for} & 1 < s
   \end{array} \right.
\label{eq:zetaBS}
\end{equation} 
The scaling length $\sigma$ is now determined by fulfilling $s = 1$ for 
$k = k_s = \pi/h$:
\begin{equation}
 \sigma \frac{\pi}{h} = 1
 \quad \Leftrightarrow \quad
 \sigma = \frac{h}{\pi}
 \quad \Rightarrow \quad
 s = \frac{h k}{\pi}
 \quad \text{and} \quad
 \rho = \frac{\pi r}{h}
\end{equation}
Here $h$ is the discretization length, chosen so the smallest physical length 
scale of the problem is sufficiently resolved.

Considering Eqs.~\eqref{eq:zetaFT} and \eqref{eq:zetaBS}, we directly see that 
the regularization function is given by:
\begin{equation}
 \zeta(\rho) = \frac{J_{\frac{d}{2}}(\rho)}{ (2 \pi \sigma^2 \rho)^{\frac{d}{2}}} 
 \label{eq:bessel}
\end{equation} 
For the case of $d = 2$ this is identical to a regularization function 
presented in the references \cite{Leonard:1980,Winckelmans:1993}.

We may now obtain the regularized Green's function in the real domain by radial 
integration of Eq.~\eqref{eq:greens poisson 2} in 
the unbounded real domain. 
For the one-dimensional case ($d=1$), we obtain:
\begin{equation}
 G(\rho) = - \dfrac{\sigma}{\pi} \big( \text{Si}( \rho ) \, \rho + \cos(\rho) \big) + C_1
\label{eq:G1D}
\end{equation}
Here $\text{Si}( \rho ) = \int_0^{\rho} \frac{\sin( t )}{ t } \ dt$ 
is the sine integral function and $C_1$ is an integration constant, 
which we may use to define a reference value.
In this work we determine $C_1$ such that we obtain an asymptotic behavior 
towards the non-regularized Green's function for large $\rho$ whatever the value of $\sigma$:
\begin{equation}
 C_1 = \frac{1}{2} L
\quad \Leftrightarrow \quad
 G(\rho) \rightarrow -\frac{1}{2} (r - L)
 \quad \text{for} \quad \rho \rightarrow \infty
\end{equation}
Here $L$ is an arbitrary reference length which typically represents the 
largest physical length scale of the problem, by which $L \gg h$.

The gradient vector of the Green's function, which may be used to combine 
the gradient operator directly in Eq.~\eqref{eq:Fourier Poisson vec}, 
is in the one-directional case ($d=1$) given by:
\begin{equation}
 K(\rho) = -\frac{ \text{Si}( \rho ) }{ \pi }
 \label{eq:K1D}
\end{equation}

For the two-dimensional case ($d=2$) the radial 
integration yields the real domain Green's function:
\begin{equation}
 G(\rho) = - \dfrac{1}{2 \pi} \text{Bi}( \rho ) + C_2
\label{eq:G2D}
\end{equation}
where $\text{Bi}( \rho ) =  \int_0^{\rho} \frac{1 - J_0( t )}{ t } \ dt$ and 
$C_2$ is an integration constant, which we may use to define a reference value. 
The integral function $\text{Bi}$ can be expressed by 
$\text{Bi}( \rho ) = \text{Ji}_0( \rho ) + \ln \left(\frac{\rho}{2} \right) + \gamma$ 
with $\gamma$ being the Euler constant and 
$\text{Ji}_0(\rho) = \int_\rho^\infty \frac{J_0( t )}{ t } \ dt$ 
denoting the Bessel integral function which for a numerical calculation 
may be approximated efficiently by Chebyshev polynomials 
(see e.g.~\cite{Luke:1975}).
Determining the integration constant $C_2$ such that we 
obtain an asymptotic behavior towards the singular Green's function for 
large $\rho$ whatever the value of $\sigma$:
\begin{equation}
 C_2 = -\frac{1}{2 \pi} \left( \text{ln} \left( \frac{2\sigma}{L} \right) - \gamma \right)
\quad \Leftrightarrow \quad
 G(\rho) \rightarrow -\frac{1}{2 \pi} \text{ln}\left(\frac{r}{L} \right)
 \quad \text{for} \quad \rho \rightarrow \infty
\end{equation}

The gradient vector of the Green's function, which may be used to calculate 
the vector field directly in Eq.~\eqref{eq:Fourier Poisson vec}, 
is in the two-directional case ($d=2$) given by:
\begin{equation}
 \boldsymbol{K}(\rho) = -K(\rho) \, \boldsymbol{e}_r 
\quad \text{with} \quad 
 K(\rho) 
         = - \frac{1}{\sigma} \frac{d G}{d\rho}
         = \frac{1}{2\pi  \sigma}  \frac{1 - J_0(\rho)}{\rho}
 \label{eq:K2D}
\end{equation}
Here $K(\rho) \geq 0$ with $K(0) = 0$ and 
$\boldsymbol{e}_r = \boldsymbol{x}/r$ is the radial unit vector.
This result is equivalent to that presented in the references 
\cite{Leonard:1980,Winckelmans:1993}.

Using the same approach for the three-dimensional domain ($d=3$), we obtain 
the Green's function:
\begin{equation}
 G(\rho) = \dfrac{1}{2 \pi^2 \sigma} \dfrac{\text{Si}(\rho)}{\rho} + C_3
\quad \text{with} \quad 
 G(0) = \dfrac{1}{2\pi^2 \sigma}
 \label{eq:G3D}
\end{equation}
The integration constant $C_3$ is again chosen such that we obtain an 
asymptotic behavior towards the singular Green's function for 
large $\rho$ whatever the value of $\sigma$:
\begin{equation}
 C_3 = 0
\quad \Leftrightarrow \quad
 G(\rho) \rightarrow \frac{1}{4 \pi r}
 \quad \text{for} \quad \rho \rightarrow \infty
\end{equation}

The gradient vector of the Green's function 
is in the three-directional case ($d=3$) obtained as:
\begin{equation}
 \boldsymbol{K}(\rho) = - K(\rho) \, \boldsymbol{e}_r 
\quad \text{with} \quad 
 K(\rho) 
         = - \frac{1}{\sigma} \frac{d G}{d\rho}
         = \dfrac{1}{2 \pi^2 \sigma^2} \dfrac{\text{Si}(\rho)-\sin(\rho)}{\rho^2}
\label{eq:K3D}
\end{equation} 
where we also have that $K(\rho) \geq 0$ and $K(0) = 0$.

\begin{figure}
\setlength{\unitlength}{0.0500bp}
\begin{picture}(8760,11000)
\put(4680,11100)
{  
 \makebox(0,0)[c]{1-Dimensional Green's functions}
}
\put(-200,7600)
{
  \begin{picture}(4680.00,3276.00)
      \put(860,640){\makebox(0,0)[r]{\strut{}0.00}}
      \put(860,1119){\makebox(0,0)[r]{\strut{}0.40}}
      \put(860,1598){\makebox(0,0)[r]{\strut{}0.80}}
      \put(860,2077){\makebox(0,0)[r]{\strut{}1.20}}
      \put(860,2556){\makebox(0,0)[r]{\strut{}1.60}}
      \put(860,3035){\makebox(0,0)[r]{\strut{}2.00}}
      \put(980,440){\makebox(0,0){\strut{}0.00}}
      \put(1537,440){\makebox(0,0){\strut{}0.50}}
      \put(2093,440){\makebox(0,0){\strut{}1.00}}
      \put(2650,440){\makebox(0,0){\strut{}1.50}}
      \put(3206,440){\makebox(0,0){\strut{}2.00}}
      \put(3763,440){\makebox(0,0){\strut{}2.50}}
      \put(4319,440){\makebox(0,0){\strut{}3.00}}
      \put(250,1837){\rotatebox{-270}{\makebox(0,0){\strut{}$G$}}}
      \put(2649,140){\makebox(0,0){\strut{}$r/h$}}
      \put(0,0){\includegraphics{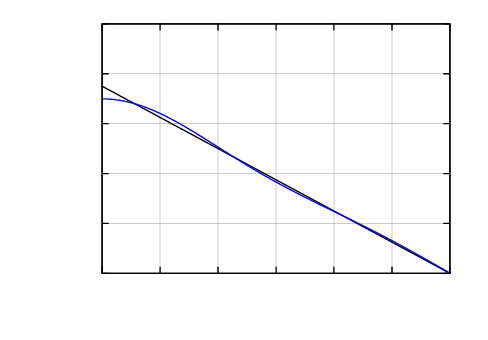}}
  \end{picture} 
}
\put(4680,7600)
{  
  \begin{picture}(4680.00,3276.00)
      \put(860,640){\makebox(0,0)[r]{\strut{}0.00}}
      \put(860,1039){\makebox(0,0)[r]{\strut{}0.10}}
      \put(860,1438){\makebox(0,0)[r]{\strut{}0.20}}
      \put(860,1838){\makebox(0,0)[r]{\strut{}0.30}}
      \put(860,2237){\makebox(0,0)[r]{\strut{}0.40}}
      \put(860,2636){\makebox(0,0)[r]{\strut{}0.50}}
      \put(860,3035){\makebox(0,0)[r]{\strut{}0.60}}
      \put(980,440){\makebox(0,0){\strut{}0.00}}
      \put(1537,440){\makebox(0,0){\strut{}1.00}}
      \put(2093,440){\makebox(0,0){\strut{}2.00}}
      \put(2650,440){\makebox(0,0){\strut{}3.00}}
      \put(3206,440){\makebox(0,0){\strut{}4.00}}
      \put(3763,440){\makebox(0,0){\strut{}5.00}}
      \put(4319,440){\makebox(0,0){\strut{}6.00}}
      \put(250,1837){\rotatebox{-270}{\makebox(0,0){\strut{}$K$}}}
      \put(2649,140){\makebox(0,0){\strut{}$r/h$}}
    \put(0,0){\includegraphics{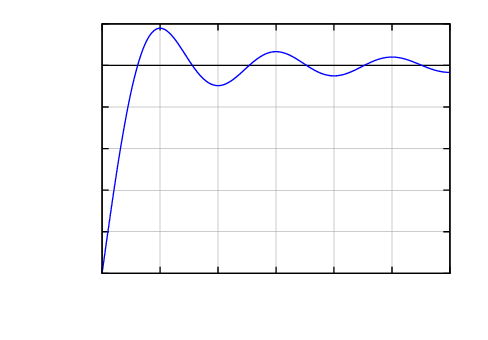}}
  \end{picture}
}

\put(4680,7300)
{  
 \makebox(0,0)[c]{2-Dimensional Green's functions}
}

\put(-200,3800)
{
  \begin{picture}(4680.00,3276.00)
      \put(860,640){\makebox(0,0)[r]{\strut{}0.00}}
      \put(860,1239){\makebox(0,0)[r]{\strut{}0.10}}
      \put(860,1838){\makebox(0,0)[r]{\strut{}0.20}}
      \put(860,2436){\makebox(0,0)[r]{\strut{}0.30}}
      \put(860,3035){\makebox(0,0)[r]{\strut{}0.40}}
      \put(980,440){\makebox(0,0){\strut{}0.00}}
      \put(1537,440){\makebox(0,0){\strut{}0.50}}
      \put(2093,440){\makebox(0,0){\strut{}1.00}}
      \put(2650,440){\makebox(0,0){\strut{}1.50}}
      \put(3206,440){\makebox(0,0){\strut{}2.00}}
      \put(3763,440){\makebox(0,0){\strut{}2.50}}
      \put(4319,440){\makebox(0,0){\strut{}3.00}}
      \put(250,1837){\rotatebox{-270}{\makebox(0,0){\strut{}$G$}}}
      \put(2649,140){\makebox(0,0){\strut{}$r/h$}}
    \put(0,0){\includegraphics{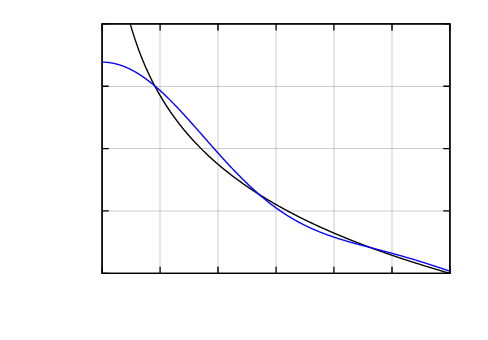}}
  \end{picture}
}
\put(4680,3800)
{  
  \begin{picture}(4680.00,3276.00)
      \put(860,640){\makebox(0,0)[r]{\strut{}0.00}}
      \put(860,1119){\makebox(0,0)[r]{\strut{}0.05}}
      \put(860,1598){\makebox(0,0)[r]{\strut{}0.10}}
      \put(860,2077){\makebox(0,0)[r]{\strut{}0.15}}
      \put(860,2556){\makebox(0,0)[r]{\strut{}0.20}}
      \put(860,3035){\makebox(0,0)[r]{\strut{}0.25}}
      \put(980,440){\makebox(0,0){\strut{}0.00}}
      \put(1537,440){\makebox(0,0){\strut{}1.00}}
      \put(2093,440){\makebox(0,0){\strut{}2.00}}
      \put(2650,440){\makebox(0,0){\strut{}3.00}}
      \put(3206,440){\makebox(0,0){\strut{}4.00}}
      \put(3763,440){\makebox(0,0){\strut{}5.00}}
      \put(4319,440){\makebox(0,0){\strut{}6.00}}
      \put(250,1837){\rotatebox{-270}{\makebox(0,0){\strut{}$K$}}}
      \put(2649,140){\makebox(0,0){\strut{}$r/h$}}
    \put(0,0){\includegraphics{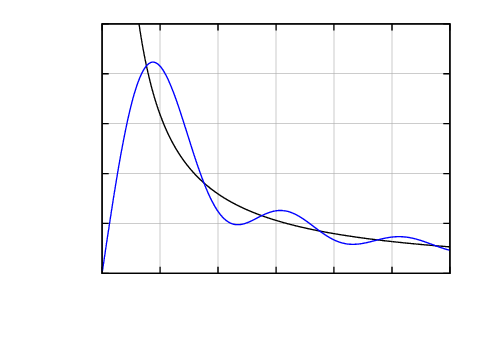}}
  \end{picture}
}

\put(4680,3500)
{  
 \makebox(0,0)[c]{3-Dimensional Green's functions}
}
\put(-200,0)
{
  \begin{picture}(4680.00,3276.00)
      \put(860,640){\makebox(0,0)[r]{\strut{}0.00}}
      \put(860,1119){\makebox(0,0)[r]{\strut{}0.05}}
      \put(860,1598){\makebox(0,0)[r]{\strut{}0.10}}
      \put(860,2077){\makebox(0,0)[r]{\strut{}0.15}}
      \put(860,2556){\makebox(0,0)[r]{\strut{}0.20}}
      \put(860,3035){\makebox(0,0)[r]{\strut{}0.25}}
      \put(980,440){\makebox(0,0){\strut{}0.00}}
      \put(1537,440){\makebox(0,0){\strut{}0.50}}
      \put(2093,440){\makebox(0,0){\strut{}1.00}}
      \put(2650,440){\makebox(0,0){\strut{}1.50}}
      \put(3206,440){\makebox(0,0){\strut{}2.00}}
      \put(3763,440){\makebox(0,0){\strut{}2.50}}
      \put(4319,440){\makebox(0,0){\strut{}3.00}}
      \put(250,1837){\rotatebox{-270}{\makebox(0,0){\strut{}$G$}}}
      \put(2649,140){\makebox(0,0){\strut{}$r/h$}}
    \put(0,0){\includegraphics{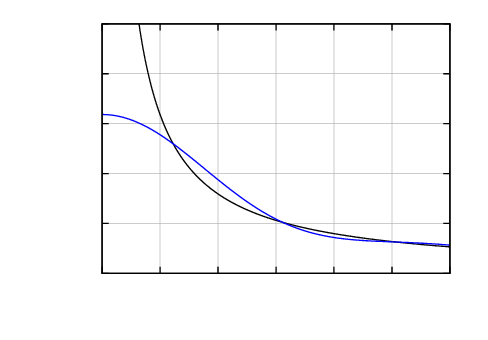}}
  \end{picture}%
}
\put(4680,0)
{  
  \begin{picture}(4680.00,3276.00)
      \put(860,640){\makebox(0,0)[r]{\strut{}0.00}}
      \put(860,1039){\makebox(0,0)[r]{\strut{}0.02}}
      \put(860,1438){\makebox(0,0)[r]{\strut{}0.04}}
      \put(860,1838){\makebox(0,0)[r]{\strut{}0.06}}
      \put(860,2237){\makebox(0,0)[r]{\strut{}0.08}}
      \put(860,2636){\makebox(0,0)[r]{\strut{}0.10}}
      \put(860,3035){\makebox(0,0)[r]{\strut{}0.12}}
      \put(980,440){\makebox(0,0){\strut{}0.00}}
      \put(1537,440){\makebox(0,0){\strut{}1.00}}
      \put(2093,440){\makebox(0,0){\strut{}2.00}}
      \put(2650,440){\makebox(0,0){\strut{}3.00}}
      \put(3206,440){\makebox(0,0){\strut{}4.00}}
      \put(3763,440){\makebox(0,0){\strut{}5.00}}
      \put(4319,440){\makebox(0,0){\strut{}6.00}}
      \put(250,1837){\rotatebox{-270}{\makebox(0,0){\strut{}$K$}}}
      \put(2649,140){\makebox(0,0){\strut{}$r/h$}}
    \put(0,0){\includegraphics{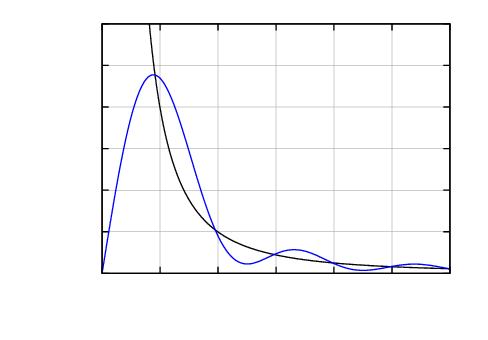}}
  \end{picture}
}
\end{picture} 

\caption{
The regularized Green's functions (blue) compared to the 
singular Green's functions (black) for the one-, two-, and three-dimensional cases. 
(Left) The Green's function $G(r)$ (Eqs.~\eqref{eq:G1D}, \eqref{eq:G2D} and \eqref{eq:G3D}). 
To determine the integration constants for the Green's functions in 
1D and 2D, an atypical value of $L = 3 h$ was here used in order to illustrate 
and emphasize the essential characteristics of the non-singular Green's functions---in 
reality, one would of course have $L \gg h$.
(Right) the radial component of the Green's function gradient $K(r)$ 
(Eqs.~\eqref{eq:K1D}, \eqref{eq:K2D} and \eqref{eq:K3D}). 
}
 \label{fig:Greens}
\end{figure}

The derived non-singular Green's functions are compared to the singular 
Green's functions in Fig.~\ref{fig:Greens}. The oscillatory behavior of the 
non-singular Green's functions is analogous to the sinc function approximation 
of the Dirac delta function, and ensures the conservation of the field moments.

For validation we have implementing the non-singular Green's functions into the 
numerical methodology used in the Poisson solver presented in references 
\cite{Hejlesen:2013,Hejlesen:2015d,Hejlesen:2016,Spietz:2018}.
Using this method we are able to obtain an accuracy in the solution limited 
only by the machine precision, when tested on the same benchmark cases as is 
used in \cite{Hejlesen:2013} (not shown). It is here emphasized that this order of 
precision is only feasible given a sufficiently smooth and resolved
right-hand-side field of the Poisson equation. 
For less smooth and resolved fields, the 
presented non-singular Green's functions still provides the highest possible 
accuracy allowed by the quadrature error. 
For the actual numerical implementation that we have used for these tests, 
the reader is referred to the provided open source software \cite{Hejlesen:2018}.

%
%

\section{Conclusion}

We have derived the non-singular Green's functions for the unbounded Poisson 
equation, in one, two and in three dimensions, using a Fourier analysis. 
By applying a radial cut-off function in the Fourier domain, 
the regularized Green's function is derived analytically in unbounded real 
domains subject to a minimum resolved length scale (e.g. a mesh size $h$). 
As a sharp cut-off function is used, smoothing errors are avoided and the 
obtained regularized Green's function achieves an optimal accuracy when used 
to solve the Poisson equation on discretized fields. 

The gradient vector of the Green's function was furthermore derived.
The obtained Green's function for the two-dimensional case corresponds 
to the one stated in references \cite{Leonard:1980,Winckelmans:1993}. 
In the present work, we formally showed the derivation of this function 
by Fourier analysis, and we extended it to the one and three-dimensional cases.

\section*{Acknowledgement}

This research did not receive any specific grant from funding agencies in the 
public, commercial, or not-for-profit sectors.
We would like to acknowledge the helpful support of Henrik Juul Spietz in 
developing the numerical methodology used in the referenced mesh based Poisson 
solver, and Denis-Gabriel Caprace for pointing out some needed improvements 
to the original manuscript.

\section*{References}
\bibliographystyle{elsarticle-num}

\begin{thebibliography}{10}
\expandafter\ifx\csname url\endcsname\relax
  \def\url#1{\texttt{#1}}\fi
\expandafter\ifx\csname urlprefix\endcsname\relax\def\urlprefix{URL }\fi
\expandafter\ifx\csname href\endcsname\relax
  \def\href#1#2{#2} \def\path#1{#1}\fi

\bibitem{Hejlesen:2019}
M.~M. Hejlesen, G. Winckelmans and J.~H. Walther. 
Non-singular Green's functions for the unbounded Poisson equation in 
one, two and three dimensions, 
Applied Mathematics Letters 89 (2019) 28--34.

\bibitem{Leonard:1980}
A.~Leonard, 
Vortex methods for flow simulation, 
J.~Comput.\ Phys. 37 (1980) 289--335.

\bibitem{Beale:1985}
J.~T.~Beale, A.~Majda, 
High order accurate vortex methods with explicit velocity kernels, 
J.~Comput.\ Phys. 58 (1985) 188--208.

\bibitem{Winckelmans:1993}
G.~S.~Winckelmans, A.~Leonard, 
Contribution to vortex particle methods for the computation of three-dimensional 
incompressible unsteady flows, 
J.~Comput.\ Phys. 109 (1993) 247--273.

\bibitem{Vico:2016}
F.~Vico, L.~Greengard, M.~Ferrando. 
Fast convolution with free-space Green’s functions. 
J. Comput. Phys. 323 (2016), 191--203.

\bibitem{Carrier:1988}
J.~Carrier, L.~Greengard, V.~Rokhlin, 
A fast adaptive multipole algorithm for particle simulations, 
SIAM J.\ Sci.\ Stat.\ Comput. 9~(4) (1988) 669--686.

\bibitem{Hockney:1970}
R.~W.~Hockney, 
The Potential calculation and some applications, 
Methods Comput.\ Phys. 9 (1970) 136--210.

\bibitem{Barnes:1986}
J.~Barnes, P.~Hut, 
A hierarchical ${O}({N}\log{N})$ force-calculation algorithm, 
Nature 324~(4) (1986) 446--449.

\bibitem{Hejlesen:2013}
M.~M.~Hejlesen, J.~T.~Rasmussen, P.~Chatelain, J.~H.~Walther, 
A high order solver for the unbounded {P}oisson equation, 
J.~Comput.\ Phys. 252 (2013) 458--467.

\bibitem{Hejlesen:2015d}
M.~M.~Hejlesen, J.~T.~Rasmussen, P.~Chatelain, J.~H.~Walther, 
High order {P}oisson solver for unbounded flows, 
Procedia IUTAM 18 (2015) 56--65.

\bibitem{Hejlesen:2016}
M.~M.~Hejlesen, 
A high order regularisation method for solving the {P}oisson
equation and selected applications using vortex methods, 
Ph.D. thesis, Technical University of Denmark (February 2016).

\bibitem{Chatelain:2010}
P.~Chatelain, P.~Koumoutsakos, 
A Fourier-based elliptic solver for vortical flows with periodic and unbounded 
directions. 
J. Comput. Phys. 229 (2010), 2425--2431.

\bibitem{Spietz:2018}
H.~J.~Spietz, M.~M.~Hejlesen, J.~H.~Walther, 
A regularization method for solving the {P}oisson equation for mixed 
unbounded-periodic domains. 
J. Comput. Phys. 356 (2018), 439--447.

\bibitem{Hejlesen:2018}
M.~M.~Hejlesen, Open source software available at: 
\texttt{https://github.com/mmhej/poissonsolver.git}

\bibitem{Rasmussen:2011c}
J.~T. Rasmussen, 
Particle methods in bluff body aerodynamics, 
{Ph.D}.\ thesis, Technical University of Denmark (October 2011).

\bibitem{Stein:1971}
E.~M. Stein, G. Weiss, Introduction to Fourier Analysis on Euclidean
Spaces, Princeton University Press, Princeton, N.J., 1971.

\bibitem{Olver:2010}
F.~W.~J. Olver, D.~W. Lozier, R.~F. Boisvert, C.~W. Clark, 
{NIST} handbook of mathematical functions, 
Cambridge University Press, New York, 2010.

\bibitem{Luke:1975}
Y.~L.~Luke, 
Mathematical functions and their approximations, 
Academic Press Inc., 1975.

\end{thebibliography}

\end{document}